\documentclass[11pt,reqno]{amsart}
\usepackage{amssymb,amsmath}\newcommand{\be}{\begin{eqnarray}}
\newcommand{\ee}{\end{eqnarray}}\newcommand{\e}{{\varepsilon}}\newcommand{\Th}{{\theta}}\newcommand{\R}{{\mathbb R}}\newcommand{\Z}{{\mathbb Z}}\newcommand{\Nat}{{\mathbb N}}\newcommand{\Hau}{{\mathcal H}}\newcommand{\E}{{\bf E\,}}\newcommand{\Cant}{{\mathcal C}}\newcommand{\lsim}{\lesssim}\newcommand{\gsim}{\gtrsim}\newcommand{\K}{{\mathcal K}}\newcommand{\M}{{\mathcal M}}\newcommand{\Pk}{{\mathcal P}}\newcommand{\Rk}{{\mathcal R}}\newcommand{\Ik}{{\mathcal I}}\newcommand{\Jk}{{\mathcal J}}\newcommand{\Kk}{{\mathcal K}}\newcommand{\Uk}{{\mathcal U}}\newcommand{\ccc}{a}
\newcommand{\la}{\lambda}\newcommand{\Proj}{\operatorname{Proj}}\newcommand{\dist}{\operatorname{dist}}\newcommand{\ci}[1]{_{{}_{\scriptstyle{#1}}}}\newcommand{\Fav}{\operatorname{Fav}}\newcommand{\supp}{\operatorname{supp}}\newtheorem{theorem}{Theorem}\newtheorem*{claim1} {Claim 1}\newtheorem*{claim2} {Claim 2}\theoremstyle{definition}\theoremstyle{remark}\numberwithin{equation}{section}\input epsf.sty
\begin{document}\thispagestyle{empty}

\title[Buffon needle probability of the four-corner Cantor set]{{The power law for the Buffon needle probability of the four-corner Cantor set}}
\author{Fedor Nazarov}\address{Fedor Nazarov, Department of  Mathematics, University of Wisconsin.
\newline{\tt nazarov@math.wisc.edu}}
\author{Yuval Peres}\address{Yuval Peres, Microsoft Research, Redmond  and Departments of Statistics and Mathematics, University of California, Berkeley.
\newline{\tt peres@microsoft.com}}
\author{Alexander Volberg}\address{Alexander Volberg, Department of  Mathematics, Michigan State University and the University of Edinburgh
{\tt volberg@math.msu.edu}\,\,and\,\,{\tt a.volberg@ed.ac.uk}}

\thanks{Research of the authors was supported in part by NSF grants  DMS-0501067 (Nazarov and Volberg) and DMS-0605166 (Peres)}
\begin{abstract}Let $\Cant_n$ be the $n$-th generation in  the construction of the middle-half Cantor set. The Cartesian square $\K_n$ of $\Cant_n$ consists of $4^n$ squares of side-length $4^{-n}$. The chance that a long needle thrown at random in the unit square will meet $\K_n$ is essentially the average length of the projections of $\K_n$, also known as the Favard length of $\K_n$.
A classical theorem of Besicovitch implies that the Favard length of  $\K_n$ tends to zero. It is still an open problem to determine its exact rate of decay. Until recently, the only explicit upper bound was $\exp(- c\log_* n)$, due to Peres and Solomyak.  ($\log_* n$ is the number of times one needs to take log to obtain a number less than $1$ starting from $n$). We obtain a power law bound by combining analytic and combinatorial ideas.\end{abstract}\maketitle

\section{{\bf Introduction}} \label{sec:intro}
\begin{figure}[htbp]\centerline{\epsfxsize=2.in \epsffile{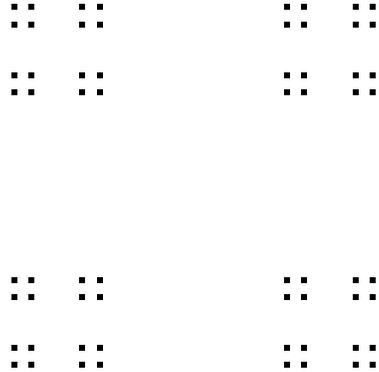}} \caption{$\K_3$, the third stage of the construction of $\K$.}\end{figure}
The four-corner Cantor set $\K$ is constructed by replacing the unit square by four sub-squares of side length $1/4$ at its corners, and iterating this operation in a self-similar manner in each sub-square.  More formally,
consider  the set  $\Cant_n$ that is the union of $2^n$ segments:

\[
\Cant_n = \bigcup_{\ a_j \in \lbrace 0,3  \rbrace, j=1,..,n}\Bigl[\sum_{j=1}^n a_j 4^{-j}, \sum_{j=1}^n a_j 4^{-j} +4^{-n}\Bigr]\,,
\]
and  let the middle half Cantor set be
\[
\Cant :=\bigcap_{n=1}^{\infty} \Cant_n\,.
\]
It can also be written as
$\Cant = \lbrace \sum_{n=1}^\infty a_n 4^{-n}:\ a_n \in \lbrace 0,3 \rbrace \rbrace.$
The four corner Cantor set $\K$ is the Cartesian square $\Cant \times \Cant$. 

Since the one-dimensional Hausdorff measure of $\K$ satisfies  $0<\Hau^1(\K)<\infty$ and the projections of $\K$ in two distinct directions have zero length, a theorem of Besicovitch (see\cite[Theorem 6.13]{falc1})  yields that the projection of $\K$ to almost every line through the origin has zero length. This is equivalent to saying that the Favard length of $\K$ equals zero.
Recall (see \cite[p.357]{besi}) that the {\bf Favard length} of a planar set $E$ is defined by\begin{equation} \label{favdef}\Fav(E) = \frac{1}{\pi} \int_0^\pi |\Proj \Rk_\Th E|\,d\Th,\end{equation}where $\Proj$ denotes the orthogonal projection from $\R^2$ to the horizontal axis,  $\Rk_\Th$ is the counterclockwise rotation by angle $\Th$, and $|A|$ denotes the Lebesgue measure of a measurable set $A \subset \R$. The Favard length of a set $E$ in the unit square has a probabilistic interpretation: up to a constant factor, it is the probability that the ``Buffon's needle,'' a long line segment dropped at random, hits $E$  (more precisely, suppose the needle's length is infinite, pick its direction uniformly at random, and then locate the needle in a uniformly chosen position in that direction, at distance at most $\sqrt{2}$ from the center of the unit square).

The set $\K_n=\Cant_n^2$  is  a union of $4^n$ squares with side length $4^{-n}$ (see Figure 1 for a picture of $\K_3$). By the dominated convergence theorem, $\Fav(\K)=0$ implies  $\lim_{n\to\infty} \Fav(\K_n)=0$. We are interested in  good estimates for $\Fav(\K_n)$ as $n\to \infty$. A lower bound $\Fav(\K_n)\ge \frac{c}{n}$ for some $c>0$ follows from  Mattila \cite[1.4]{mattila1}.   Peres and Solomyak \cite{PS} proved that

\[
\Fav(\K_n) \le C \exp[-\ccc\log_*n]\ \ \ \ \mbox{for all}\ \ n \in \Nat,\] where  \[\label{def-logstar}\log_* n = \min\left\lbrace k\ge 0:\ \underbrace{\log\log\ldots \log}_{k}n \le 1 \right \rbrace\,.
\]

This result can be viewed as an attempt to make  a quantitative statement out of a qualitative Besicovitch projection theorem \cite{besi}, \cite{T}, using this canonical example of the Besicovitch irregular set.

It is very interesting to see what are quantitative analogs of Besicovitch theorem in general. The reader can find more of that in \cite{T}.

We now state our main result, which improves this upper bound to a power law.

\begin{theorem}
\label{th-corner}
For every $\delta>0$, there exists $C>0$ such that
\[\Fav(\K_n) \le C n^{\delta-1/6}\ \ \ \ \mbox{for all}\ \ n\in \Nat.\]
\end{theorem}

\noindent {\bf Remarks.}
\begin{itemize}
\item The $1/6$ in the exponent is certainly not optimal, and, indeed, can be  improved slightly with the methods of this paper. However, a bound decaying faster than $O\Bigl(n^{-1/4}\Bigr)$ would require new ideas.
\item In \cite{PS}, Theorem 2.2,  a  random analog of the Cantor set $\K$ is analyzed, and it is shown that, with high probability,   the Favard length of the $n$-th stage in the construction has upper and lower bounds  that are constant multiples of $n^{-1}$. However, it is not clear to us whether $\Fav(\K_n)$ also decays at this rate.
\item It follows from the results of Kenyon \cite{kenyon} and Lagarias and Wang \cite{lawang} that $|\Proj \Rk_\Th\K|=0$ for all $\Th$ such that $\tan\Th$ is irrational.
As noted in \cite{PS}, this information does not seem to help obtain an upper bound for  $\Fav(\K_n)$.
\item The set $\K$ was one of the first examples of sets of positive length and zero analytic capacity, see \cite{david} for a survey. The asymptotic behavior of the analytic capacity of $\K_n$ was determined in 2003 by Mateu, Tolsa and Verdera~\cite{MTV}, it is equivalent to $\frac{1}{\sqrt{n}}$.
\end{itemize}

\medskip

It will be convenient to translate $\K_n$ so that its convex hull is the unit square centered at the origin. Due to the symmetries of the square, one can average over $\theta \in (0,\frac \pi 4)$ in the definition  (\ref{favdef}) of $\Fav(\K \ci n)$. After translation, the projection of  $\Rk \ci{\Th} \Bigl(\K \ci n-(\frac{1}{2},\frac{1}{2})\Bigr)$ to  the horizontal axis is the union of $4^n$ intervals of length $4^{-n}(\cos \theta + \sin \theta)$ centered at the points $\sum_{k=0}^{n-1}4^{-k}\xi_k$, where

\[\xi_k\in \Big \lbrace \pm \frac{3 \sqrt{2}}{8} \cos (\frac \pi 4 - \theta),\pm \frac{3 \sqrt{2}}{8} \sin (\frac \pi 4 - \theta) \Big \rbrace.\]

Let now  $t=\tan(\frac \pi 4-\theta)\in[0,1]$. Since $\frac {\sqrt{2}} 2 \leq \cos(\frac \pi 4 - \theta) \leq 1$ on $(0, \frac \pi 4)$, the length of the projection $\Proj \Rk_\Th (\K \ci n)$ is comparable to the length of the union of $4^n$ intervals of length $4^{-n} \rho$ centered at the points $\sum_{k=0}^{n-1}4^{-k}\xi_k$ with $\xi_k \in \lbrace \pm 1, \pm t  \rbrace$, where $\rho = \rho(\Th) = \frac{8}{3\sqrt{2}} (1 + \tan(\frac \pi 4 - \Th))$. The exact value of $\rho(\Th)$ is of no importance, the only thing that matters is that it is separated from both 0 and $+\infty$. We shall also need the function $f_n$ that is the product of $\frac{1}{\rho}$ and the sum of the characteristic functions of these intervals. In other words,

\[f_n = \nu^{(n)} * \frac {4^n} {\rho} \chi\ci{[-\frac{\rho}{2} 4^{-n}, \frac{\rho}{2}4^{-n}],}\] where\[\nu^{(n)}=*_{k=0}^{n-1}\nu_k, \quad \text{and} \quad \nu_k=\frac 14[\delta_{-4^{-k}}+\delta_{-4^{-k}t}+\delta_{4^{-k}t}+\delta_{4^{-k}}]\,.
\]

Geometrically,  $f_n$ is  (up to minor rescaling) the number of squares whose projections contain a given point. Finally, since $\lvert \frac{dt}{d\theta} \rvert = \frac{1}{\cos^2(\frac \pi 4 - \theta)}$ is between 1 and 2 for all $\theta \in [0, \frac{\pi} {4})$, we can replace averaging over $\theta$ with that over $t$.

\section{Fourier-analytic part.}

In what follows, we will use $\asymp$ and $\lsim$, $\gsim$ to denote, respectively, equality or the corresponding inequality up to some positive multiplicative constant. Let $K, S$ be large positive numbers. Our first aim is to show that there exists a power $p>0$ (we'll see that any $p > 4$ fits) such that the measure of the set\[E=\left\lbrace t\in[0,1]:\max_{1\le n\le (KS)^p}\int_{\R} f_n^2\le K\right \rbrace\] is at most $\frac 1S$. Suppose not.
Let $N$ be the least even integer exceeding $\frac12 (KS)^p$. For every $t\in E$, we must have\[K\ge \int_{\mathbb R}|f\ci N(x)|^2\,dx \asymp \int_{\mathbb R}|\widehat f\ci N(y)|^2\,dy\gsim \int_{1}^{4^{N/2}}|\widehat{\nu}^{(N)}(y)|^2\,dy\,,\] because $\psi=\frac{4^N}{\rho}  \chi\ci{\left[-\frac{\rho}{2} 4^{-N} \, , \,\frac{\rho}{2}4^{-N}\right]}$ satisfies $\widehat{\psi}(y) \gsim 1$
for all $\vert y \vert < 4^{N/2}$ if $N$ is sufficiently large. Thus

\[\frac 1{|E|}\int_E \Bigl[\sum_{n=1}^{N/2}\int_{4^{n-1}}^{4^{n}}|\widehat{\nu}^{(N)}(y)|^2\,dy\Bigr]\,dt\le K
\]
and for each $m\le N/2$, there exists $n\le N/2$ satisfying

\[\frac 1{|E|}\int_E\Bigl[\int_{4^{n-m}}^{4^{n}}|\widehat \nu^{(N)}(y)|^2\,dy\Bigr]\,dt\le \frac {4Km}N\,.
\]

Thus\[E_*=\Bigl\{ t \in E \, : \, \int_{4^{n-m}}^{4^{n}}|\widehat{\nu}^{(N)}(y)|^2\,dy\le \frac {8Km}{N} \Bigr\}\] satisfies $|E_*| \ge |E|/2$. Our assumption on $E$ implies that  $|E_*| \ge \frac{1}{2S}$. Now for $y \in [4^{n-m},4^n]$, we have\[\lvert \widehat{\nu}^{(N)}(y)\rvert^2\asymp \prod_{k=0}^{n}\Bigl\lvert \frac{\cos 4^{-k}y+\cos 4^{-k}ty}{2}\Bigr\rvert^2\,,\] since the remaining terms (that correspond to $k \in [n+1,N]$) in the product converge geometrically to 1. Making the change of variable $y\mapsto 4^ny$, we get

\[\int_{4^{n-m}}^{4^{n}}\lvert \widehat \nu^{(N)}(y)\rvert ^2\,dy\asymp4^n\int_{4^{-m}}^1 \Bigl\lvert\prod_{k=0}^{n}\frac{\cos 4^{k}y+\cos 4^{k}ty}{2}\Bigr\rvert^2\,dy\,.
\]

Now split the last product into
$$P_1(y)=\prod_{k=0}^{m}\frac{\cos 4^{k}y+\cos 4^{k}ty}{2} \; \;  \mbox{ \rm and} \; \;P_2(y)=\prod_{k=m+1}^{n}\frac{\cos 4^{k}y+\cos 4^{k}ty}{2} \,.$$ Consider the integral

\[\int_{4^{-m}}^1 |P_2(y)|^2\,dy\]
first. Writing the cosines as  sums of exponentials,  we have

\[P_2(y)=4^{m-n}\sum_{j=1}^{4^{n-m}}e^{i\la_j y} \,,\]
where $\{\lambda_j\}_{j=1}^{4^{n-m}}$ are the sums of all subsets of $\{\pm 4^k, \pm 4^k t \, : \, k \in [m+1,n]\}$. For $t \in E_* \subset E$, the definition of $E$ yields that

\[\int_{\R}\Bigl(\sum_j\chi\ci{[\la_j - \frac \rho2 4^m, \la_j + \frac \rho2 4^m]}\Bigr)^2\le K\cdot 4^n
\]
(this is equivalent to $\int_{\R}f_{n-m}^2\le K$).
The last inequality can be viewed as a separation condition on the spectrum, so one can hope that a variation of Salem's trick should allow us to conclude that

\[\int_{4^{-m}}^1|P_2(y)|^2\,dy\gsim 4^{m-n}\,
,\]
provided that $L=4^m$ is chosen appropriately. We shall choose $m$ such that $4^m=L$ is a large constant multiple of $K$. Since  $|P_2(-y)|=|P_2(y)|$, we can integrate over $[-1,1]\setminus[-L^{-1},L^{-1}]$. Consider the function $g$ given by

\[g(y) = (1 - |y|)_{+} - 2 (1 - L^{-1})(1 - \tfrac L2 |y|)_+ +(1 -2 L^{-1})(1 - L|y|)_+.\]

\begin{figure}
\centerline{\epsfxsize=2.5in \epsffile{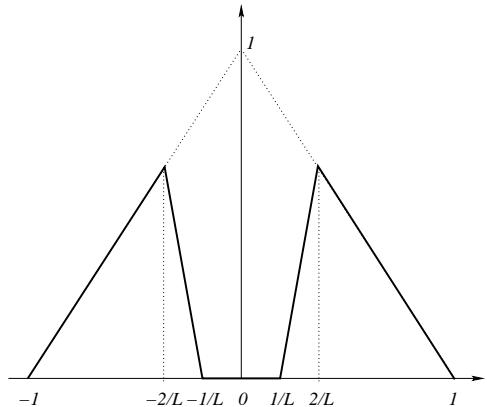}} \caption{Triangle kernel function.}
\end{figure}

Note that  $g$ is even, $0 \leq g \leq 1$,  $\supp g\subset [-1,1]\setminus[-L^{-1},L^{-1}]$ and  $\int_{-1}^1 g\ge \frac 12$ if $L$ is not too small.  Now, let $h$ denote  ``the triangle function" that is $1$ at $0$, vanishes on $\R\setminus (-1,1)$ and is linear on $[-1,0]$ and on $[0,1]$.
Then
\[
g(y) = h(y) - 2(1-L^{-1})h(\tfrac{L}2 y)  + (1-2L^{-1})h(Ly),.
\]
As $\widehat h (\la) = 2\frac{1-\cos \la}{\la^2}\in [0, \frac{C}{1 +\lambda^2}]$, we get
\[
\widehat g (\la) \ge -\frac{C}{L} \cdot \frac{1}{1 + (\lambda/L)^2}\,.
\]
So we got the estimate
\[
\widehat g (\la) \geq -C\frac{L}{\la^2 + L^2}
\]
with some numerical constant $C$.


\medskip

Denote $M=4^{n-m}$. Let us call $k\in\lbrace 1,\dots, M \rbrace$ {\em good\/} if\[C\sum_j\frac{L}{L^2+(\la_j-\la_k)^2}\le\frac 18\,.\] Then

\begin{multline*}\int_{[-1,1]\setminus[-L^{-1},L^{-1}]}\Bigl|\sum_{k}e^{i\la_k y}\Bigr|^2\,dy\,\ge \int_{\R}g(y)\Bigl|\sum_{k}e^{i\la_k y}\Bigr|^2\,dy\,\\\ge \sum_{ \lbrace k:k\text{ is good}\rbrace }\frac 12+\int_{\R}g(y)\Bigl|\sum_{ \lbrace k: k\text{ is bad}\rbrace}e^{i\la_k y}\Bigr|^2\,dy\,-{2\sum_{ \lbrace k: k\text{ is good} \rbrace} C\sum_j\frac{L}{L^2+(\la_j-\la_k)^2}}\\\ge \frac 14\#\lbrace k:k\text{ is good} \rbrace\,.
\end{multline*}

Now we need only to show that the number of good indices is comparable to $M$.
To this end, note that we have the condition

\[\int_{\R}\Bigl(\sum_j\chi\ci{[\la_j - \frac \rho2 L \, , \, \la_j+  \frac \rho2 L]}(\la)\Bigr)^2\,d\la\le MLK\,.\]

Convolving with the Poisson kernel $\Pk\ci L(\la)=\frac 1\pi\frac {L}{L^2+\la^2}$ and taking into account that\[\chi\ci{[\la_j - \frac \rho2 L,\la_j+  \frac \rho2 L]}*\Pk\ci L\ge c L\Pk\ci L( \cdot - \la_j )\] with $c > 0$ (here we use that $\rho$ stays bounded away from 0 and $+\infty$), we get

\[L^2\int_{\R}\Bigl[\sum_j \Pk\ci L(\la-\la_j)\Bigr]^2\,d\la\le C'MLK\,,\] but
\[\int_{\R}\Pk\ci L(\la-\la_j)\Pk\ci L(\la-\la_k)\,d\la\ge c'\Pk\ci L(\la_j-\la_k)\,.\]

Thus\[c'\sum_{j,k}\Pk\ci L(\la_j-\la_k)\le C'MKL^{-1}\,\] and\[\#\lbrace k:k\text{ is bad} \rbrace \le \frac{8C\pi C'}{c'}(KL^{-1})M \le \frac M2,\] provided that $L \ge \frac {16 C\pi C'}{c'} K$.
Therefore, indeed, $\int_{4^{-m}}^1 |P_2(y)|^2\,dy\, \gsim 4^{m-n}$.

\medskip

The danger is that this large integral can be completely killed when the integrand is multiplied by $|P_1|^2$. Note that

\[\frac{\cos 4^k y+\cos 4^k ty}2=\cos 2^{-1}4^k(y+ty)\cos 2^{-1}4^k(y-ty)\] so\[P_1(y)=\prod_{k=0}^m \cos 2^{-1}4^k(y+ty)\cos 2^{-1}4^k(y-ty).\] Using the formula\[2\cdot 4^m \sin(\frac u2)\prod_{\ell=0}^{2m}\cos 2^{\ell-1}u=\sin 4^m u\,.\]
we conclude that

\[|P_1(y)|\gsim 4^{-2m}|\sin 4^m(y+ty)|\cdot|\sin 4^m(y-ty)|\,.\]

This can be small only if $\sin 4^m(y+ty)$ or $\sin 4^m(y-ty)$ is small.
For $\delta \in (0,1)$, denote by $\Ik \ci \delta$ the union of intervals  of length $4^{-m}\delta$ centered at the points $\frac{\pi \ell}{4^m}$, $\ell \in \Z$.   Define  $\omega(t; \delta)$ by

\[\omega(t; \delta) = \lbrace  y \in (4^{-m},1): \; y+ty \in \Ik \ci \delta \;  \mbox{ \rm or} \; y - ty \in \Ik \ci \delta \rbrace.\]

We would like to estimate $\int_{\omega(t; \delta)} |P_2(y)|^2\,dy$ from above. This may be a hard task for an individual $t\in E_*$, but we can bound the average fairly easily.
We have
\begin{multline*}
\frac 1{|E_*|}\int_{E_*}\Bigl(\int_{\omega(t; \delta,)}|P_2(y)|^2\,dy\Bigr)\,dt\le 2S \int_0^1\Bigl(\int_{\omega(t; \delta,)}|P_2(y)|^2\,dy\Bigr)\,dt\\\lsim 2S\int_{[4^{-m},1] \cap \Ik \ci \delta}\Bigl( \prod_{k=m+1}^n\cos^2 2^{-1}4^k u\Bigr)\,\frac{du}{u+v}\int_{[0,1] }\Bigl( \prod_{k=m+1}^n\cos^2 2^{-1}4^k v\Bigr)\,dv\\+ 2S\int_{[4^{-m},1] }\Bigl( \prod_{k=m+1}^n\cos^2 2^{-1}4^k u\Bigr)\,\frac{du}{u+v}\int_{[0,1]\cap \Ik \ci \delta }\Bigl( \prod_{k=m+1}^n\cos^2 2^{-1}4^k v\Bigr)\,dv
\end{multline*}
where $u = y+ty$ and $v = y-ty$.

Using the formula $\cos^2\alpha=\frac 12(1+\cos 2\alpha)$ and the inequality $\frac1{u+v}\le L\,du\,dv$, we can estimate the last expression by
$$
\E:=C\,S\,L\cdot 4^{m-n}\Bigl[\int_\mathcal{I} \prod_{k=m+1}^n (1+\cos 4^k u) \,du\Bigr]\cdot\Bigl[\int_{[0,1]}\prod_{k=m+1}^n (1+\cos 4^k v) \,dv\Bigr] \,.
$$

Above we observed that\[\prod_{k=m+1}^n (\cos^2 2^{-1} 4^k u) = 2^{m-n}\prod_{k=m+1}^n (1+\cos 4^k u) =: 2^{m-n} R(u)\,.\]

We want to see now that

$$
\E\le c\,S\,L\cdot 4^{m-n}\sqrt\delta\,.
$$

To this end we notice that the Riesz product $R(u)$ is $\frac{\pi}{4^m}\,$-periodic.

 Note also that, for any interval $\Jk$ of length $4^{-m} \frac{\pi} {4^j}$ ($j \in \Z_+$), we have

  \[\int_\Jk R(u)\,du = \int_{\Jk} R_1(u)R_2(u)\,du,
  \]
where $R_1(u) = \prod_{k=m+1}^{m+j} (1+\cos 4^k u)$ and  $R_2(u) = \prod_{k=m+j+1}^{n} (1+\cos 4^k u)$. Observe that $R_1(u) \leq 2^j$ for all $u$ and $R_2(u)$ is $\frac{\pi}{4^{m+j}}\,$-periodic, so
  \[\int_\Jk R_2(u)\,du = \frac 1 {4^{m+j}} \int_{0}^\pi R_2(u)\,du = \frac{\pi}{4^{m+j}}.
  \]
  Thus  $\int_{\Jk}R(u)\,du \leq \frac{\pi}{2^j }4^{-m}$.
  Choose $j$ in such a way that $\delta \asymp 4^{-j}$.

  It follows that, for each constituting interval $\Jk$ of $\Ik \ci \delta$, we have $\int_{\Jk}R(u)\,du \lsim 4^{-m} \sqrt{\delta}$.

  \[\int_{[4^{-m},1] \cap \Ik \ci \delta} R(u)\,du \lsim 4^m\cdot 4^{-m} \sqrt{\delta} \lsim \sqrt{\delta}\,.
  \]
  In conjunction with the  estimate $\int_{[0,1] \cap \Ik \ci \eta} R(v)\,du \lsim 1$, we finally get

  $$
\E\le c\,S\,L\cdot 4^{m-n}\sqrt\delta\,.
$$

The resulting estimate is much less than $4^{m-n}$ if $\delta$ is much less than $S^{-2}L^{-2}$. Thus, for at least one $t\in E$, we must have
(recall that $L=4^m$)

 $$
 \int_{[L^{-1},1]\setminus\Omega(t)}|P_2(y)|^2\,dy\ge c 4^{m-n}
 $$
 and, thereby, (if we remember that $K$ was a small constant times $4^m$)

 $$ \int_{[L^{-1},1]}|P_1(y)|^2|P_2(y)|^2\,dy\ge 4^{-4m}(S^{-2}L^{-2})^4\cdot 4^{m-n}\ge c S^{-8} K^{-11}4^{-n}\,. $$

  Thus, if $p>12$ then our choice of $N$ at the beginning of the proof gives $N> (KS)^{12+\e}/2> (KS)^{12+\e/2}$,  hence $\dfrac {2K\log K}N$ is much less than $S^{-8}K^{-11}$, and we get a contradiction.

  However, we promised to show that  $N> (KS)^{4+\e}$ already leads to  a contradiction. To do this, we make our considerations more elaborate, but we follow the same lines. In fact, let us consider

  \[\Omega(t; \delta,\eta) = \lbrace  y \in (4^{-m},1): \; y+ty \in \Ik \ci \delta \;  \mbox{ \rm and} \; y - ty \in \Ik \ci \eta \rbrace.\]

We changed the word ``or" in the definition of $\omega(t;\delta)$ by the word ``and" in the definition of $\Omega(t;\delta,\eta)$. This will allow us to make a subtler estimate. Notice that
$$
\{y: |\sin 4^m(y+ty)|\cdot|\sin 4^m(y-ty)|\le 2^{-l}\} \subset \bigcup_{k=0}^\ell \Omega (t; 2^{-k}, 2^{-\ell+k+1})\,.
$$

We would like to estimate $\int_{\Omega(t; \delta,\eta)} |P_2(y)|^2\,dy$ from above. As before
we have
\begin{multline*}
\frac 1{|E_*|}\int_{E_*}\Bigl(\int_{\omega(t; \delta,)}|P_2(y)|^2\,dy\Bigr)\,dt\le 2S \int_0^1\Bigl(\int_{\Omega(t; \delta, \eta)}|P_2(y)|^2\,dy\Bigr)\,dt\\\lsim 2S\int_{[4^{-m},1] \cap \Ik \ci \delta}\Bigl( \prod_{k=m+1}^n\cos^2 2^{-1}4^k u\Bigr)\,\frac{du}{u}\int_{[0,1]\cap \Ik \ci \eta }\Bigl( \prod_{k=m+1}^n\cos^2 2^{-1}4^k v\Bigr)\,dv\end{multline*}
where $u = y+ty$ and $v = y-ty$ as before.

We already introduced $R(u)=\prod_{k=m+1}^n (1+\cos 4^k u)$ and established the following estimate

$$
\int_{[0,1]\cap \mathcal{I}_{\eta}} R(v)\,dv \lsim \sqrt\eta\,.
$$

Now we can estimate
$$
\int_{[4^{-m},1]\cap \mathcal{I}_{\delta}} R(u)\frac{du}{u} \lsim \sum_{1\le j\le \frac1{\pi}4^m}
\frac{4^m}{\pi j}\cdot 4^{-m}\sqrt\delta \le \sqrt\delta\,m\,.
$$

Therefore, we obtain

  \[\frac 1{|E_*|}\int_{E_*}\Bigl(\int_{\Omega(t; \delta, \eta)}|P_2(y)|^2\,dy\Bigr)\,dt\lsim Sm \sqrt{\delta \eta}\,4^{m-n}\,.
  \]

  Let us denote $\Omega\ci\ell (t):=\bigcup_{k=0}^\ell \Omega (t; 2^{-k}, 2^{-\ell+k+1})$.  We know now that

  \[\frac 1{|E_*|}\int_{E_*}\Bigl(\int_{\Omega \ci \ell (t)}|P_2(y)|^2\,dy\Bigr)\,dt\lsim Sm \ell \cdot 2^{-\ell/2}\,4^{m-n}.
  \]
  If $Sm \ell \cdot 2^{-\ell/2}$ is a small constant (much less than 1), then it follows   that there exists $t \in E_*$ such that

  \[
  \int_{[4^{-m},1] \setminus \Omega \ci \ell (t)}|P_2(y)|^2\,dy \ge c\cdot 4^{m-n}.
  \]
  But $\Omega\ci\ell(t)$ contains
  $
\{y: |\sin 4^m(y+ty)|\cdot|\sin 4^m(y-ty)|\le 2^{-l}\} $. This means that
$|P_1| \gsim 4^{-2m} 2^{-\ell}$ on $(4^{-m},1) \setminus \Omega \ci \ell$, so for this $t$, we have
  \[
  \int_{4^{-m}}^{1} |P_2(y)|^2\,dy \geq 4^{-4m} 2^{-2 \ell} 4^{m-n}.
  \]
  If $\frac{Km}N$ is much less than $4^{-4m} 2^{-2 \ell} 4^{m}$, we get a contradiction.
  Since $4^m \asymp K$, we see that we get a contradiction if it is possible to find $\ell$ such that $Sm \ell \cdot 2^{-\ell/2}$ is much less than 1 and $N$ is much greater than $K^4 m 2^{2\ell}$ simultaneously. If $N > (KS)^{4+ \gamma}$ with $\gamma > 0$, we can take $2^{\ell/2} \asymp (SK)^{\frac \gamma 8}S$, thus finishing the proof of our claim with any $p>4$.

  \section{Combinatorial part.}

  Fix the rotation angle $\theta$ and  some large positive integer $N$. As before, let $F_n(x)$ be the number of the squares in $\Rk \ci \theta \K \ci n$ whose projections to the horizontal axis contain $x$.
  Define
  \[F_{*}(x) = \max_{0 \leq n \leq N}  F_n(x).
  \]
  Our  key observation is the following inequality:  for any positive integers $K,M$, we have
  \[\mu\lbrace F_*\ge 4KM \rbrace\le 108 K \mu\lbrace F_*\ge K \rbrace\mu\lbrace F_*\ge M \rbrace,
  \]
  where $\mu$ denotes the usual Lebesgue measure on the real line.

  \vspace{5mm}

  \noindent{ \bf{Proof:}}
  For each point $x \in \R$ where $F_*(x)\ge 2K$, choose the least $n= n(x)$ for which $F_n(x)\ge 2K$ and mark all the squares in $\Rk \ci \theta \K\ci n$ whose projections contain $x$. Note that the number of such squares for a given point $x$ cannot exceed $4K$: otherwise we would have $F_{n-1}(x) \geq 2K$, which contradicts our choice of $n$.
  Now unmark all marked squares that are contained in larger marked squares and consider the family of the remaining maximal marked squares. The desired inequality is immediately implied by the following two claims:

  \begin{claim1}In order to reach the level $4KM$ at $x$, we have to reach level $M$ in at least one maximal marked square whose projection contains $x$.
  \end{claim1}

  \begin{claim2}The sum of side lengths of all maximal marked squares does not exceed $108 K\mu \lbrace F_* \geq K  \rbrace$.
  \end{claim2}

  \vspace{5mm}

  \noindent{ \bf{Proof of Claim 1:}}
  Obviously, in order to reach the level $4KM$ at $x$, one has to reach the level $M$ in at least one of the squares of generation $n(x)$ whose projection contains $x$ (recall that there are not more than $4K$ such squares!) Each such square is contained in some maximal marked square, whence the claim.

  \vspace{5mm}

  \noindent{ \bf{Proof of Claim 2:}}
  Consider all 4-adic intervals $I \subset \R$ such that $I$ intersects a projection of some maximal marked square $Q$ whose side length is at least $|I|$. Clearly, the union of all such intervals contains the projections of all maximal marked squares. Now consider maximal intervals $I$ with this property. Clearly, each such maximal interval $I$ intersects the projection of some maximal marked square $Q$ with side length $|I|$, but no projection of a maximal marked square with a larger side length.

  Now let us estimate the sum of side lengths of the maximal marked squares intersecting one such maximal 4-adic interval $I$. Let $\sigma = \sin \theta + \cos \theta$. Note that each maximal marked square whose projection intersects $I$ is contained in some square of generation $\log_4 \frac 1{|I|}$ with side length $|I|$, whose projection intersects $I$. Since the projection of each such square is contained in $(2\sigma + 1)I$, having more than $\frac{2\sigma +1}\sigma \cdot 4K$ such squares would imply existence of a point $x \in I$ that is contained in more than $4K$ projections of squares of generation $\log_4 \frac 1{|I|}$.
  But this implies that there are at least $2K$ squares of the previous generation above $x$, so $n(x) \leq \log_4 \frac 1{|I|}-1$ and there exists a marked square of side length greater than $|I|$ whose projection intersects $I$. The maximal marked square containing it has at least the same side length and its projection still intersects $I$.  But this contradicts maximality of $I$.

  Since the maximal marked squares are disjoint, the sum of side lengths of maximal marked squares contained in one square of generation $\log_4 \frac 1{|I|}$ does not exceed $|I|$. Hence the sum of side lengths of all maximal marked squares whose projections intersect $I$ is at most $\frac{2\sigma +1}\sigma \cdot 4K|I| \leq 12 K|I|$. Thus, the total sum of side lengths of all maximal marked squares is at most $12K \sum_{I} |I| = 12K |U_{I} I|$, because the maximal intervals are disjoint.

  Now let $I$ be one of our maximal intervals and let $Q_1$ be a maximal marked square with side length $|I|$, whose projection intersects $I$. Since $Q_1$ is a marked square, there exists a point $x$ and $2K-1$ other squares $Q_2, \ldots, Q_{2K}$ of side length $|I|$ such that $\Proj Q_j \ni x$ for all $j = 1, \ldots, 2K$. Choosing $K$ such squares whose centers lie on one side of $x$, we see that there exists an interval $J$ of length $\sigma/2$ containing $x$. such that $F_* \geq F \ci{\log_4 \frac 1{|I|}} \geq K$ on $J$. Since $\dist(x,I) \leq \sigma |I|$, we have $I \subset \frac{5\sigma + 4}{\sigma} J \subset 9J$. Hence, if $J'$ is the constituting interval of the set $\lbrace F_* \geq K  \rbrace$, containing $J$, we also have $I \subset 9J'$. Therefore, $|\bigcup_I I| \leq 9\mu \lbrace F_* \geq K \rbrace$ and we are done.

  Now fix $\theta$, $K$, and $N$. Let $\nu = \mu \lbrace F_* \geq K  \rbrace$.
  By induction, we get\[\mu \lbrace F_* \geq (4K)^j K  \rbrace \leq [108 K \nu]^j \nu, \quad j = 1,2,\ldots..\]
  Hence, for all $n=0,1,\ldots,N$, we get
  \begin{multline*}
  \int_ \R f_n^2 = \int_{\lbrace f_n \leq K \rbrace } f^2_n + \int_{\lbrace K\leq f_n \leq 4K^2 \rbrace} f^2_n+ \sum\limits_{j \geq 1} \int_{\lbrace (4K)^j K \leq f_n \leq (4K)^{j+1}K \rbrace } f^2_n\\\leq \sqrt2 K + 16 K^4 \nu + \sum_{j \geq 1}[108K\nu]^j (4K)^{2j} 16 K^4 \nu \leq 2K,
  \end{multline*}
  provided that $108\cdot16K^3 \nu \leq \frac 12$. The Fourier-analytic part implies that the measure of all angles $\theta$ with this property is less than some absolute constant times $\frac K{N^{1/4 - \gamma}}$ with arbitrarily small $\gamma > 0$.

  Assume now that $\nu > 32^{-1}\cdot 108^{-1}\cdot K^{-3}$. For each point $x \in \R$ where $F_*(x) \geq K$, choose some $n \in \lbrace 0,1,\ldots,N \rbrace$ for which $F_n(x) \geq K$ and mark all squares of the $n$-th generation whose projections contain $x$. Now, in the $N$-th generation, color green all squares contained in the marked squares. Let $\varphi$ be the sum of the characteristic functions of projections of green squares and let $\Xi$ be he union of the projections of all marked squares.
  We want to show first that
  \[\Xi \subset \lbrace y \in \R: \M \varphi(y) \geq \frac{K}{4}  \rbrace,
  \]
  where $\M\varphi(y) = \sup_{r>0} \frac1{2r} \int_{y-r}^{y+r} \varphi(s)\, ds$ is the central Hardy-Littlewood maximal function.
  Indeed, if $y \in \Xi$, then the vertical line through $y$ intersects at least one marked square $Q_1$. Thus, there exists $x \in \R$ and $K-1$ other marked squares $Q_2, \ldots, Q_K$ of the same size as $Q_1$, such that $x \in \Proj Q_j$ for all $j = 1, \ldots, K$. Now, the interval $J$ centered at $y$ of length $4| \Proj Q_1|$ contains all the projections of the squares $Q_1, \ldots, Q_K$. Integral $\int_{J} \varphi$ is then not less than the sum of all lengths of the projections of the green squares contained in $Q_1, \ldots, Q_k$, which is $K| \Proj Q_1|$.
  Hence, $\M\varphi(y) \geq \frac 1{|J|} \int_{J} \varphi \geq \frac K4$. Using the weak type $L^1$ estimate for the Hardy-Littlewood maximal function, we conclude that
  \[\mu (\Xi) \lsim \frac{1}{K} \int_{\R} \varphi.
  \]
  Since $F_*(x) \geq K$ implies $x \in \Xi$, we deduce that $\int_\R\varphi \gsim K\nu \gsim K^{-2}$, i.e., there are at least $cK^{-2}4^{N}$ green squares.

  On the other hand, $\Xi$ contains the projections of all green squares, and $\int_\R \varphi$  is (up to a constant factor)  the sum of all side lengths of all green squares. Thus,  the length of the projection of the union of all green squares is at most $\frac{C}{K}$ times the sum of their side lengths.
  The net outcome of the previous construction is that in the $N$-th generation rotated Cantor square $\Rk \ci \theta \Kk \ci N$, we have $\Uk \gsim K^{-2} 4^N$ green squares, whose projections overlap a lot (more precisely, their total projection is only about $\frac{1}{K}$ times their total side length) and $4^N - \Uk$ other (white) squares about which we know nothing.
  This just gives the estimate $\frac{C}{K} \,\Uk \cdot 4^{-N} + \sqrt 2 \cdot 4^{-N}(4^N - \Uk)$ for the total length of the projection of $\Rk \ci \theta \Kk \ci N$, which doesn't look very impressive. But here is where the self-similarity comes into play.

  Let us repeat the construction of green squares in each of the white squares (this will bring us to the consideration of $K_{2N}$ instead of $K_N$). Now we will have $\Uk \cdot 4^N$ squares contained in the original green squares, which still give us the projections $\leq \frac{C}{K} \Uk \cdot 4^{-N}$, but we shall also have $(4^N - \Uk)\Uk$ new small green squares and the total length of their projection will be ${\leq \frac{C}{K} \Uk (4^N - \Uk)\cdot 4^{-2N}}$. Thus the total length of the projections of all these squares will be at most ${ \frac{C}{K} \Uk \cdot 4^{-N}[1 + (1 - \frac{\Uk}{4^N})]}$.
  To this we should add $\sqrt{2}\cdot 4^{-2N}(4^N - \Uk)^2 = \sqrt{2} (1 - \frac{\Uk}{4^N})^2$, which is the trivial upper bound for the total projection of the remaining $(4^N - \Uk)^2$ squares in $\Rk \ci \theta \Kk \ci {2N}$. Proceeding to $ \Kk \ci {3N}$ in a similar manner, we shall get\[\lvert \Proj \Rk \ci \theta \Kk \ci {3N}\rvert \leq\frac{C}{K}\, \Uk \cdot 4^{-N}[1 + (1 - \frac{\Uk}{4^N})+(1 - \frac{\Uk}{4^N})^2 ]+\sqrt{2} (1 - \frac{\Uk}{4^N})^3,\]
  and so on.
  By the time we reach $\Rk \ci \theta \Kk \ci {XN}$ with a large positive integer $X$, we shall get\[\lvert \Proj \Rk \ci \theta \Kk \ci {XN}\rvert \leq\frac{C}{K} \frac{\Uk} {4^{N}} \sum_{\ell = 0} ^{X-1} (1 - \frac{\Uk}{4^N})^\ell+\sqrt{2} (1 - \frac{\Uk}{4^N})^X.\]
  The first term does not exceed $\frac{C}{K} \frac{\Uk} {4^{N}} \sum_{\ell = 0} ^{\infty} (1 - \frac{\Uk}{4^N})^\ell = \frac{C}{K}$, while the second is at most $\sqrt{2} e^{- 4^{-N} \Uk X}$, which is less than $\frac{\sqrt {2}}{K}$ if $4^{-N} \Uk X > \log K$, i.e., if $X$ is much greater than $K^2\log K$.
  The moral of the story is that, given two positive integers $K$ and $N$, we can find an exceptional set of measure  $\lsim \frac{K}{N^{1/4 -\gamma}}$, such that for all $\theta$ outside this set, we have $\lvert \Proj \Rk \ci \theta \Kk \ci {XN}\rvert \leq \frac{1}{K}$ for all integers $X$ that are much greater than $K^2 \log K$.

  The last result can be restated as follows: If $K$, $S$ are large enough and {$N \geq K^p S^q$} with $p>6$, $q>4$, then
  \[\mu \Bigl \lbrace  \theta \in (0, \frac{\pi}4): \lvert \Proj \Rk \ci \theta \Kk \ci {N}\rvert \geq \frac{C}{K} \Bigl \rbrace \lsim \frac{1}{S}.
  \]
  This gives us the weak type inequality
  \[\mu \Big \lbrace  \theta \in (0, \frac{\pi}4): \lvert \Proj \Rk \ci \theta \Kk \ci {N}\rvert \geq t \Big \rbrace \lsim (N^{-1} t^{-p})^{1/q},
  \]
  provided that $N^{-1} t^{-p}$ is much less than 1.
  Combining it with the trivial estimate $\mu \lbrace  \theta \in (0, \frac{\pi}4): \ldots \rbrace \leq \frac{\pi}4$ for all other $t$, we finally get:
  \begin{multline*}\int_0^{\frac{\pi}4} \lvert \Proj \Rk \ci \theta \Kk \ci {N}\rvert \,d\theta = \int_0^{\infty} \mu\Bigl \lbrace \theta \in (0, \frac{\pi}4): \lvert \Proj \Rk \ci \theta \Kk \ci {N}\rvert
  \geq t \Bigr \rbrace\, dt \,\\\lsim \int_0^{C N^{-1/p}} 1 \,dt\, + \int_{C N^{-1/p}}^{\infty} N^{-1/q} t^{-p/q}\,dt\,\\\lsim N^{-\frac{1}{p}} + N^{-\frac{1}{q}} N^{\frac{1}{p}(\frac{p}{q}-1)} = 2 N^{-\frac{1}{p}},
  \end{multline*}
  finishing the proof.

  \section{$h$-Hausdorff measures of the projections.}

  If  a function $h$ is increasing, continuous and $h(0)=0$, we can define Hausdorff measure $\Hau_h$ on compact set by the usual procedure. When $h(t)=t$, this is exactly the Hausdorff measure $\Hau^1$ of dimension $1$. We know that $\Hau^1$ measure of almost all projections of $1/4$ Cantor set is zero, and the Hausdorff dimension of almost every projection is $1$. We can get more information about these projections by measuring their $\Hau_h$ using a more refined scale of gauge functions than just powers of $t$.  Namely, the main result obtained in this paper readily implies the following corollary. Consider the gauge function $h(t) = t(\log \frac1t)^c$ with small positive $c$. We have proved

  \begin{theorem}
  If $c$ is sufficiently small ($c\in (0, 1/6)$) then almost every projection of the four corner Cantor set $\K$ has zero $\Hau_h$ measure.
  \end{theorem}

  \bibliographystyle{amsplain}
\end{document}